# QUANTITATIVE BOUNDS ON CONVERGENCE OF TIME-INHOMOGENEOUS MARKOV CHAINS

By R. Douc, E. Moulines and Jeffrey S. Rosenthal[1]

*École Polytechnique, École Nationale Supérieure des Télécommunications and University of Toronto*

Convergence rates of Markov chains have been widely studied in recent years. In particular, quantitative bounds on convergence rates have been studied in various forms by Meyn and Tweedie [*Ann. Appl. Probab.* **4** (1994) 981–1101], Rosenthal [*J. Amer. Statist. Assoc.* **90** (1995) 558–566], Roberts and Tweedie [*Stochastic Process. Appl.* **80** (1999) 211–229], Jones and Hobert [*Statist. Sci.* **16** (2001) 312–334] and Fort [Ph.D. thesis (2001) Univ. Paris VI]. In this paper, we extend a result of Rosenthal [*J. Amer. Statist. Assoc.* **90** (1995) 558–566] that concerns quantitative convergence rates for time-homogeneous Markov chains. Our extension allows us to consider $f$-total variation distance (instead of total variation) and time-inhomogeneous Markov chains. We apply our results to simulated annealing.

## 1. Time-homogeneous case.

1.1. *Introduction.* Let $P$ be a Markov transition kernel defined on a general state space $(\mathcal{X}, \mathcal{B}(\mathcal{X}))$. Denote by $P^k$ the corresponding $k$-step transition kernel. For $\xi$ a probability measure on $\mathcal{B}(\mathcal{X})$ and $f$ a Borel function, define $\xi P(A) = \int \xi(dy) P(y, A)$ and $Pf(x) = \int P(x, dy) f(y)$.

For $f: \mathcal{X} \to [1, \infty)$, the *$f$-total variation* or *$f$-norm* of a signed measure $\mu$ on $\mathcal{B}(\mathcal{X})$ is defined as

$$\|\mu\|_f := \sup_{|\phi| \leq f} |\mu(\phi)|.$$

When $f \equiv 1$, the $f$-norm is the total variation norm, which is denoted $\|\mu\|_{\mathrm{TV}}$. Our goal is to find explicit bounds on rates of convergence of $\xi P^n - \xi' P^n$

Received April 2002; revised September 2003.
[1]Supported in part by NSERC of Canada.

to zero. In the special case in which $P$ has a stationary distribution $\pi$, this corresponds to bounding the convergence of $\xi P^n$ to $\pi$. Our results extend and sharpen the nonquantitative results developed, for example, by Meyn and Tweedie [(1993), Chapters 15 and 16], where one typically finds conditions under which there exists some *rate function* $r(n)$ such that $r(n)\|P^n(x,\cdot) - \pi\|_f \to 0$ as $n \to \infty$.

The problem of getting explicit bounds on $\|P^n(x,\cdot) - \pi\|_f$ has received much attention in recent years, motivated by control of convergence for Markov chain Monte Carlo and operation research problems [see, e.g., Jones and Hobert (2001)]. Most of the results available cover only total variation bound [see Rosenthal (1995) and Roberts and Tweedie (1999)]. To the best of our knowledge, the only explicit bound in $f$-total variation distance was given by Meyn and Tweedie [(1994), Theorem 2.3]. This bound is based on the Nummelin splitting construction and depends in a very intricate way on the constants of the kernel. In this section, we use a different approach, based on *coupling*. We obtain a bound (Theorem 2) which is simple, very generally applicable and, although not tight, does improve on the work of Meyn and Tweedie [(1994), Theorem 2.3].

1.2. *Assumptions and lemma.* Let $a \wedge b = \min(a,b)$ and $a \vee b = \max(a,b)$. To use the coupling construction, we first need a set where coupling may occur. We make the following assumption:

(A1) There exist a set $\bar{C} \subset \mathcal{X} \times \mathcal{X}$, a constant $\varepsilon > 0$ and a family of probability measures $\{\nu_{x,x'}, (x,x') \in \bar{C}\}$ on $\mathcal{X}$ with

(1) $\qquad P(x,A) \wedge P(x',A) \geq \varepsilon \nu_{x,x'}(A) \qquad \forall A \in \mathcal{B}(\mathcal{X}),\ (x,x') \in \bar{C}.$

Following Bickel and Ritov (2001), we call $\bar{C}$ a $(1,\varepsilon)$-*coupling set*. For simplicity, only one-step minorization is considered in this paper. Adaptations to $m$-step minorization can be carried out as in Rosenthal (1995). We note that condition (1) is in many cases satisfied by setting $\bar{C} = C \times C$, where $C$ is a so-called *pseudo-small* set. Recall that a subset $C \subset \mathcal{X}$ is $(1,\varepsilon)$-*pseudo-small* if there exist a constant $\varepsilon > 0$ and a family of probability measure $\{\nu_{x,x'}, (x,x') \in C \times C\}$ with $P(x,\cdot) \wedge P(x',\cdot) \geq \varepsilon \nu_{x,x'}(\cdot)$ for all $(x,x') \in C \times C$ [see Roberts and Rosenthal (2001)]. We stress that $C$ is a subset of $\mathcal{X}$ and that, despite the obvious similarity, a $(1,\varepsilon)$-pseudo-small set is not a $(1,\varepsilon)$-coupling set. Recall finally that a set $C$ is $(1,\varepsilon)$-small if it is $(1,\varepsilon)$-pseudo-small with the same minorizing probability measure $\nu = \nu_{x,x'}$ for all $(x,x') \in C \times C$. The primary motivation for using $(1,\varepsilon)$-coupling set is that the usual pairwise coupling argument can be used without change and that, in some cases detailed below, $(1,\varepsilon)$-coupling sets can be significantly larger than the product of $(1,\varepsilon)$-pseudo-small sets.

To introduce the coupling construction, some additional definitions are required. Let $\bar{R}$ be a Markov transition kernel that satisfies, for all $(x, x') \in \bar{C}$ and all $A \in \mathcal{B}(\mathcal{X})$,

(2)
$$\bar{R}(x, x'; A \times \mathcal{X}) = (1 - \varepsilon)^{-1}(P(x, A) - \varepsilon \nu_{x,x'}(A)),$$
$$\bar{R}(x, x'; \mathcal{X} \times A) = (1 - \varepsilon)^{-1}(P(x', A) - \varepsilon \nu_{x,x'}(A)).$$

For example, we can set, for $(x, x') \in \bar{C}$,

$$\bar{R}(x, x'; A \times A') = ((1 - \varepsilon)^{-1}(P(x, A) - \varepsilon \nu_{x,x'}(A)))$$
$$\times ((1 - \varepsilon)^{-1}(P(x', A') - \varepsilon \nu_{x,x'}(A'))),$$

but other trickier constructions may also be considered. Similarly, let $\bar{P}$ be a Markov transition kernel on $\mathcal{X} \times \mathcal{X}$ such that, for $(x, x') \in \bar{C}$ and all $A, A' \in \mathcal{B}(\mathcal{X})$,

(3) $\quad \bar{P}(x, x'; A \times A') = (1 - \varepsilon)\bar{R}(x, x'; A \times A') + \varepsilon \nu_{x,x'}(A \cap A'),$

and satisfies, for $(x, x') \notin \bar{C}$ and all $A \in \mathcal{B}(\mathcal{X})$,

(4) $\quad \bar{P}(x, x'; A \times \mathcal{X}) = P(x, A) \quad \text{and} \quad \bar{P}(x, x'; \mathcal{X} \times A) = P(x', A).$

For example, we can once again set, for $(x, x') \notin \bar{C}$, $\bar{P}(x, x'; A \times A') = P(x, A)P(x', A')$, to get that $\bar{P}$ satisfies (4) for all $(x, x') \in \mathcal{X} \times \mathcal{X}$.

Define the product space $\mathsf{Z} = \mathcal{X} \times \mathcal{X} \times \{0, 1\}$ and the associated product sigma algebra $\mathcal{B}(\mathsf{Z})$. We define on the space $(\mathsf{Z}^{\mathbb{N}}, \mathcal{B}(\mathsf{Z})^{\otimes \mathbb{N}})$ a Markov chain $(Z_n := (X_n, X'_n, d_n), n \geq 0)$. Indeed, given $Z_n$, we construct $Z_{n+1}$ as follows. If $d_n = 1$, then draw $X_{n+1} \sim P(X_n, \cdot)$, and set $X'_{n+1} = X_{n+1}$ and $d_{n+1} = 1$. If $d_n = 0$ and $(X_n, X'_n) \in \bar{C}$, flip a coin with probability of heads $\varepsilon$. If the coin comes up heads, then draw $X$ from $\nu_{X_n, X'_n}(\cdot)$, and set $X_{n+1} = X'_{n+1} = X$ and $d_{n+1} = 1$. If the coin comes up tails, then draw $(X_{n+1}, X'_{n+1})$ from the residual kernel $\bar{R}(X_n, X'_n; \cdot)$ and set $d_{n+1} = 0$. If $d_n = 0$ and $(X_n, X'_n) \notin \bar{C}$, then draw $(X_{n+1}, X'_{n+1})$ according to the kernel $\bar{P}(X_n, X'_n; \cdot)$ and set $d_{n+1} = 0$. Here $d_n$ is called a *bell variable*; it indicates whether the chains have coupled ($d_n = 1$) or not ($d_n = 0$) by time $n$.

For $\mu$ a probability measure on $\mathcal{B}(\mathsf{Z})$, denote by $\mathbb{P}_\mu$ the probability measure induced on $(\mathsf{Z}^{\mathbb{N}}, \mathcal{B}(\mathsf{Z})^{\otimes \mathbb{N}})$ by the Markov chain $(Z_n, n \geq 0)$ with initial distribution $\mu$. The corresponding expectation operator is denoted by $\mathbb{E}_\mu$. It is then easily checked that $(X_n, n \geq 0)$ and $(X'_n, n \geq 0)$ are each marginally updated according to the transition kernel $P$; that is, for any $n$, for any initial distributions $\xi$ and $\xi'$, and for any $A, A' \in \mathcal{B}(\mathcal{X})$,

(5)
$$\mathbb{P}_{\xi \otimes \xi' \otimes \delta_0}(Z_n \in A \times \mathcal{X} \times \{0, 1\}) = \xi P^n(A),$$
$$\mathbb{P}_{\xi \otimes \xi' \otimes \delta_0}(Z_n \in \mathcal{X} \times A' \times \{0, 1\}) = \xi' P^n(A'),$$

where $\delta_x$ is the Dirac measure centered on $x$ and $\otimes$ is the tensor product of measures. Define the *coupling time* $T = \inf\{k \geq 1; d_k = 1\}$ (with the convention $\inf \varnothing = \infty$). Let $P^*$ be the Markov kernel defined, for all $(x, x') \in \mathcal{X} \times \mathcal{X}$ and all $A \in \mathcal{B}(\mathcal{X} \times \mathcal{X})$, by

$$
(6) \quad P^*(x, x'; A) = \begin{cases} \bar{P}(x, x'; A), & \text{if } (x, x') \notin \bar{C}, \\ \bar{R}(x, x'; A), & \text{if } (x, x') \in \bar{C}. \end{cases}
$$

For $\mu$ a probability measure on $\mathcal{X} \times \mathcal{X}$, denote by $\mathbb{P}_\mu^*$ and $\mathbb{E}_\mu^*$ the probability and the expectation induced by the Markov chain on $\mathcal{X} \times \mathcal{X}$ with initial distribution $\mu$ and transition kernel $P^*$.

LEMMA 1. *Assume* (A1). *Then, for any $n \geq 0$ and any nonnegative Borel function $\phi: (\mathcal{X} \times \mathcal{X})^{n+1} \to \mathbb{R}^+$, we have*

$$\mathbb{E}_{\xi \otimes \xi' \otimes \delta_0}\{\phi(\bar{X}_0, \ldots, \bar{X}_n)\mathbf{1}(d_n = 0)\} = \mathbb{E}_{\xi \otimes \xi'}^*\{\phi(\bar{X}_0, \ldots, \bar{X}_n)(1-\varepsilon)^{N_{n-1}}\},$$

*where $\bar{X}_i := (X_i, X_i')$, $N_i := \sum_{j=0}^i \mathbf{1}_{\bar{C}}(\bar{X}_j)$ and $N_{-1} := 0$.*

PROOF. We first verify that the result holds for all functions $\phi(\bar{x}_0, \ldots, \bar{x}_n) = \prod_{i=0}^n \psi_i(\bar{x}_i)$, where $\bar{x}_i := (x_i, x_i')$ and $(\psi_i, i \geq 0)$ are nonnegative Borel functions on $\mathcal{B}(\mathcal{X} \times \mathcal{X})$. The proof is by induction. For $n = 0$, the result is obvious. Assume that the result holds up to order $n-1$ for some $n \geq 1$. We have

$$\mathbb{E}_{\xi \otimes \xi' \otimes \delta_0}\{\phi(\bar{X}_0, \ldots, \bar{X}_n)\mathbf{1}(d_n = 0)\}$$

$$= \mathbb{E}_{\xi \otimes \xi' \otimes \delta_0}\left\{\prod_{i=0}^{n-1} \psi_i(\bar{X}_i)\mathbf{1}_{\bar{C}^c}(\bar{X}_{n-1})\psi_n(\bar{X}_n)\mathbf{1}(d_n = 0)\right\}$$

$$+ \mathbb{E}_{\xi \otimes \xi' \otimes \delta_0}\left\{\prod_{i=0}^{n-1} \psi_i(\bar{X}_i)\mathbf{1}_{\bar{C}}(\bar{X}_{n-1})\psi_n(\bar{X}_n)\mathbf{1}(d_n = 0)\right\},$$

where $\bar{C}^c := \mathcal{X} \setminus \bar{C}$. Define $\mathcal{G}_k = \sigma(Z_i = (\bar{X}_i, d_i), 0 \leq i \leq k)$. Note that, for $n \geq 1$,

$$\mathbb{E}\{\psi_n(\bar{X}_n)\mathbf{1}(d_n = 0)|\mathcal{G}_{n-1}\}\mathbf{1}_{\bar{C}^c}(\bar{X}_{n-1})\mathbf{1}(d_{n-1} = 0)$$
$$= \bar{P}\psi_n(\bar{X}_{n-1})\mathbf{1}_{\bar{C}^c}(\bar{X}_{n-1})\mathbf{1}(d_{n-1} = 0).$$

Since $N_{n-2}\mathbf{1}_{\bar{C}^c}(\bar{X}_{n-1}) = N_{n-1}\mathbf{1}_{\bar{C}^c}(\bar{X}_{n-1})$ and $\bar{P}(x, x'; \cdot) = P^*(x, x'; \cdot)$ for $(x, x') \notin \bar{C}$, we have, under the induction assumption,

$$\mathbb{E}_{\xi \otimes \xi' \otimes \delta_0}\left\{\prod_{i=0}^{n-1} \psi_i(\bar{X}_i)\mathbf{1}_{\bar{C}^c}(\bar{X}_{n-1})\psi_n(\bar{X}_n)\mathbf{1}(d_n = 0)\right\}$$

$$= \mathbb{E}_{\xi \otimes \xi' \otimes \delta_0}\left\{\prod_{i=0}^{n-1} \psi_i(\bar{X}_i)\mathbf{1}_{\bar{C}^c}(\bar{X}_{n-1})\bar{P}\psi_n(\bar{X}_{n-1})\mathbf{1}(d_{n-1} = 0)\right\}$$

(7)
$$= \mathbb{E}^*_{\xi \otimes \xi'}\left\{\prod_{i=0}^{n-1} \psi_i(\bar{X}_i)\mathbf{1}_{\bar{C}^c}(\bar{X}_{n-1})P^*\psi_n(\bar{X}_{n-1})(1-\varepsilon)^{N_{n-1}}\right\}$$
$$= \mathbb{E}^*_{\xi \otimes \xi'}\left\{\prod_{i=0}^{n} \psi_i(\bar{X}_i)\mathbf{1}_{\bar{C}^c}(\bar{X}_{n-1})(1-\varepsilon)^{N_{n-1}}\right\}.$$

Similarly, note that

$$\mathbb{E}\{\mathbf{1}(d_n=0)\psi_n(\bar{X}_n)|\mathcal{G}_{n-1}\}\mathbf{1}_{\bar{C}}(\bar{X}_{n-1})\mathbf{1}(d_{n-1}=0)$$
$$= (1-\varepsilon)\bar{R}\psi_n(\bar{X}_{n-1})\mathbf{1}(d_{n-1}=0).$$

Since $(N_{n-2}+1)\mathbf{1}_{\bar{C}}(\bar{X}_{n-1}) = N_{n-1}\mathbf{1}_{\bar{C}}(\bar{X}_{n-1})$ and $\bar{R}(x,x';\cdot) = P^*(x,x';\cdot)$ for all $(x,x') \in \bar{C}$, the induction assumption implies

$$\mathbb{E}_{\xi \otimes \xi' \otimes \delta_0}\left\{\prod_{i=0}^{n-1} \psi_i(\bar{X}_i)\mathbf{1}_{\bar{C}}(\bar{X}_{n-1})\psi_n(\bar{X}_n)\mathbf{1}(d_n=0)\right\}$$
$$= (1-\varepsilon)\mathbb{E}_{\xi \otimes \xi' \otimes \delta_0}\left\{\prod_{i=0}^{n-1} \psi_i(\bar{X}_i)\mathbf{1}_{\bar{C}}(\bar{X}_{n-1})\bar{R}\psi_n(\bar{X}_{n-1})\mathbf{1}(d_{n-1}=0)\right\}$$

(8)
$$= \mathbb{E}^*_{\xi \otimes \xi'}\left\{\prod_{i=0}^{n-1} \psi_i(\bar{X}_i)\mathbf{1}_{\bar{C}}(\bar{X}_{n-1})P^*\psi_n(\bar{X}_{n-1})(1-\varepsilon)^{N_{n-1}}\right\}$$
$$= \mathbb{E}^*_{\xi \otimes \xi'}\left\{\prod_{i=0}^{n} \psi_i(\bar{X}_i)\mathbf{1}_{\bar{C}}(\bar{X}_{n-1})(1-\varepsilon)^{N_{n-1}}\right\}.$$

Thus, the two measures on $\mathcal{B}(\mathcal{X} \times \mathcal{X})^{\otimes(n+1)}$ defined, respectively, by

$$A \mapsto \mathbb{E}_{\xi \otimes \xi' \otimes \delta_0}\{\mathbf{1}_A(\bar{X}_0,\ldots,\bar{X}_n)\mathbf{1}(d_n=0)\} \quad \text{and}$$
$$A \mapsto \mathbb{E}^*_{\xi \otimes \xi'}\{\mathbf{1}_A(\bar{X}_0,\ldots,\bar{X}_n)(1-\varepsilon)^{N_{n-1}}\}$$

are equal on the monotone class $\mathcal{C} = \{A : A = A_0 \times \cdots \times A_n, A_i \in \mathcal{B}(\mathcal{X} \times \mathcal{X})\}$ and thus these two measures coincide on the product sigma algebra, which concludes the proof. □

1.3. *Main time-homogeneous result.* Let $f: \mathcal{X} \to [1,\infty]$ and let $\phi: \mathcal{X} \to \mathbb{R}$ be any Borel function such that $\sup_{x \in \mathcal{X}} |\phi(x)|/f(x) < \infty$. Using (5), the classical coupling inequality [see, e.g., Thorisson (2000), Chapter 2, Section 3] implies that

$$|\xi P^n \phi - \xi' P^n \phi| = |\mathbb{E}_{\xi \otimes \xi' \otimes \delta_0}\{\phi(X_n) - \phi(X'_n)\}|$$
$$= |\mathbb{E}_{\xi \otimes \xi' \otimes \delta_0}\{(\phi(X_n) - \phi(X'_n))\mathbf{1}(d_n=0)\}|$$
$$\leq \left(\sup_{x \in \mathcal{X}} |\phi(x)|/f(x)\right)\mathbb{E}_{\xi \otimes \xi' \otimes \delta_0}\{(f(X_n) + f(X'_n))\mathbf{1}(d_n=0)\}.$$

By Lemma 1,
$$\mathbb{E}_{\xi \otimes \xi' \otimes \delta_0}\{(f(X_n) + f(X'_n))\mathbf{1}(d_n = 0)\} = \mathbb{E}^*_{\xi \otimes \xi'}\{(f(X_n) + f(X'_n))(1-\varepsilon)^{N_{n-1}}\}.$$

Thus, the following key coupling inequality holds:

$$\begin{aligned}
&|\xi P^n \phi - \xi' P^n \phi| \\
&\qquad \leq \left(\sup_{x \in \mathcal{X}} |\phi(x)|/f(x)\right) \mathbb{E}^*_{\xi \otimes \xi'}\{(f(X_n) + f(X'_n))(1-\varepsilon)^{N_{n-1}}\}.
\end{aligned} \tag{9}$$

To bound the term on the right-hand side of (9), we need a *drift condition* outside $\bar{C}$ for the kernel $P^*$:

(A2) There exist a function $\bar{V}: \mathcal{X} \times \mathcal{X} \to [1, \infty)$ and constants $b$ and $\lambda$, $0 < \lambda < 1$, such that

$$P^* \bar{V} \leq \lambda \bar{V} + b \mathbf{1}_{\bar{C}}. \tag{10}$$

THEOREM 2. *Assume* (A1) *and* (A2). *Let* $f: \mathcal{X} \to [1, \infty)$ *be a function which satisfies* $f(x) + f(x') \leq 2\bar{V}(x, x')$ *for all* $(x, x') \in \mathcal{X} \times \mathcal{X}$. *Then, for all* $j \in \{1, \ldots, n+1\}$ *and for all initial probability measures* $\xi$ *and* $\xi'$ *on* $\mathcal{X}$,

$$\|\xi P^n - \xi' P^n\|_{\mathrm{TV}} \leq 2(1-\varepsilon)^j \mathbf{1}(j \leq n) + 2\lambda^n B^{j-1}(\xi \otimes \xi')(\bar{V}), \tag{11}$$

$$\begin{aligned}
\|\xi P^n - \xi' P^n\|_f &\leq 2(1-\varepsilon)^j (b(1-\lambda)^{-1} + \lambda^n (\xi \otimes \xi')(\bar{V}))\mathbf{1}(j \leq n) \\
&\quad + 2\lambda^n B^{j-1}(\xi \otimes \xi')(\bar{V}),
\end{aligned} \tag{12}$$

*where*

$$B = 1 \vee \left((1-\varepsilon)\lambda^{-1} \sup_{(x, x') \in \bar{C}} \bar{R}\bar{V}(x, x')\right).$$

PROOF. For any $j \in \{1, \ldots, n+1\}$, we have

$$\begin{aligned}
&\mathbb{E}^*_{\xi \otimes \xi'}\{(f(X_n) + f(X'_n))(1-\varepsilon)^{N_{n-1}}\} \\
&\qquad \leq \mathbb{E}^*_{\xi \otimes \xi'}\{(f(X_n) + f(X'_n))(1-\varepsilon)^{N_{n-1}}\mathbf{1}(N_{n-1} \geq j)\} \\
&\qquad\quad + 2\mathbb{E}^*_{\xi \otimes \xi'}\{\bar{V}(\bar{X}_n)(1-\varepsilon)^{N_{n-1}}\mathbf{1}(N_{n-1} < j)\}.
\end{aligned} \tag{13}$$

Consider the first term on the right-hand side of (13). We have

$$\begin{aligned}
&\mathbb{E}^*_{\xi \otimes \xi'}\{(f(X_n) + f(X'_n))(1-\varepsilon)^{N_{n-1}}\mathbf{1}(N_{n-1} \geq j)\} \\
&\qquad \leq (1-\varepsilon)^j \mathbb{E}^*_{\xi \otimes \xi'}\{(f(X_n) + f(X'_n))\}.
\end{aligned} \tag{14}$$

If $f \equiv 1$, then $\mathbb{E}^*_{\xi \otimes \xi'}\{(f(X_n) + f(X'_n))\} = 2$. Otherwise, by repeated application of the drift condition (A2), we have

$$(P^*)^n \bar{V} \leq \lambda (P^*)^{n-1} \bar{V} + b \leq \lambda^n \bar{V} + b \sum_{k=0}^{n-1} \lambda^k \leq \lambda^n \bar{V} + b/(1-\lambda).$$

Since $f(x) + f(x') \leq 2\bar{V}(x,x')$, we get

$$\mathbb{E}^*_{\xi\otimes\xi'}\{(f(X_n) + f(X'_n))\} \leq 2\mathbb{E}^*_{\xi\otimes\xi'}\{\bar{V}(\bar{X}_n)\} \leq 2\lambda^n(\xi\otimes\xi')(\bar{V}) + 2b/(1-\lambda).$$

Consider the second term on the right-hand side of (13). Denote for $s \geq 0$,

$$M_s := \lambda^{-s} B^{-N_{s-1}} \bar{V}(\bar{X}_s)(1-\varepsilon)^{N_{s-1}}.$$

We show that $(M_s, s \geq 0)$ is an $(\mathcal{F}, \mathbb{P}^*_{\xi\otimes\xi'})$ supermartingale, where $\mathcal{F} := (\mathcal{F}_s := \sigma(\bar{X}_i, i \leq s), s \geq 0)$. The definition of $N_s$ and the drift condition (A2) imply

$$\mathbf{1}_{\bar{C}^c}(\bar{X}_s)N_s = \mathbf{1}_{\bar{C}^c}(\bar{X}_s)N_{s-1} \quad \text{and} \quad \mathbf{1}_{\bar{C}^c}(\bar{X}_s)P^*\bar{V}(\bar{X}_s) \leq \mathbf{1}_{\bar{C}^c}(\bar{X}_s)\lambda\bar{V}(\bar{X}_s).$$

Thus, we have

$$\mathbb{E}^*\{M_{s+1}|\mathcal{F}_s\}\mathbf{1}_{\bar{C}^c}(\bar{X}_s)$$
(15)
$$= \lambda^{-(s+1)} B^{-N_s} P^*\bar{V}(\bar{X}_s)(1-\varepsilon)^{N_s}\mathbf{1}_{\bar{C}^c}(\bar{X}_s)$$
$$= \lambda^{-(s+1)} B^{-N_{s-1}} P^*\bar{V}(\bar{X}_s)(1-\varepsilon)^{N_{s-1}}\mathbf{1}_{\bar{C}^c}(\bar{X}_s) \leq M_s\mathbf{1}_{\bar{C}^c}(\bar{X}_s).$$

By definition, $\sup_{(x,x')\in\bar{C}} \bar{R}\bar{V}(x,x') \leq \lambda(1-\varepsilon)^{-1}B$. Since by construction $\mathbf{1}_{\bar{C}}P^*\bar{V} = \mathbf{1}_{\bar{C}}\bar{R}\bar{V}$, we have

$$\mathbb{E}^*\{\bar{V}(\bar{X}_{s+1})|\mathcal{F}_s\}\mathbf{1}_{\bar{C}}(\bar{X}_s) = \bar{R}\bar{V}(\bar{X}_s)\mathbf{1}_{\bar{C}}(\bar{X}_s) \leq \lambda(1-\varepsilon)^{-1}B\mathbf{1}_{\bar{C}}(\bar{X}_s).$$

Since $\mathbf{1}_{\bar{C}}(\bar{X}_s)N_s = \mathbf{1}_{\bar{C}}(\bar{X}_s)(N_{s-1} + 1)$, we have

$$\mathbb{E}^*\{M_{s+1}|\mathcal{F}_s\}\mathbf{1}_{\bar{C}}(\bar{X}_s)$$
(16)
$$\leq \lambda^{-(s+1)} B^{-1} B^{-N_{s-1}}(1-\varepsilon)^{N_{s-1}+1}\lambda(1-\varepsilon)^{-1}B\mathbf{1}_{\bar{C}}(\bar{X}_s) \leq M_s\mathbf{1}_{\bar{C}}(\bar{X}_s).$$

Equations (15) and (16) show that $(M_s, s \geq 0)$ is an $(\mathcal{F}, \mathbb{P}^*_{\xi\otimes\xi'})$ supermartingale. By the optional stopping theorem, $\mathbb{E}^*_{\xi\otimes\xi'}\{M_n\} \leq \mathbb{E}^*_{\xi\otimes\xi'}\{M_0\}$. Since $B \geq 1$, we have $\mathbf{1}(N_{n-1} < j) \leq B^{j-1}B^{-N_{n-1}}$, which implies

$$\mathbb{E}^*_{\xi\otimes\xi'}\{\bar{V}(\bar{X}_n)(1-\varepsilon)^{N_{n-1}}\mathbf{1}(N_{n-1} < j)\}$$
(17)
$$\leq \lambda^n B^{j-1}\mathbb{E}^*_{\xi\otimes\xi'}\{\lambda^{-n}B^{-N_{n-1}}\bar{V}(\bar{X}_n)(1-\varepsilon)^{N_{n-1}}\}$$
$$\leq \lambda^n B^{j-1}\mathbb{E}^*_{\xi\otimes\xi'}\{M_n\} \leq \lambda^n B^{j-1}\xi\otimes\xi'(\bar{V}).$$

By combining (14) and (17) for $f \equiv 1$, we have

$$\mathbb{E}^*_{\xi\otimes\xi'}\{(f(X_n) + f(X'_n))(1-\varepsilon)^{N_{n-1}}\}$$
$$\leq 2(1-\varepsilon)^j\mathbf{1}(j \leq n) + 2\lambda^n B^{j-1}\xi\otimes\xi'(\bar{V})$$

and (11) follows from (9). Similarly, for $f$ such that $f(x) + f(x') \leq 2\bar{V}(x,x')$, we have

$$\mathbb{E}^*_{\xi\otimes\xi'}\{(f(X_n) + f(X'_n))(1-\varepsilon)^{N_{n-1}}\}$$
$$\leq 2(1-\varepsilon)^j(\lambda^n(\xi\otimes\xi')(\bar{V}) + b/(1-\lambda)) + 2\lambda^n B^{j-1}\xi\otimes\xi'(\bar{V})$$

and (12) follows from (9). $\square$

1.4. *Application to convergence to stationarity.* If $P$ has a stationary distribution $\pi$, (i.e., if $\pi P = \pi$), then we can choose $\xi' = \pi$. Then $\pi P^n = \pi$ for all $n$ and, hence, the results (11) and (12) allow us to bound $\|\xi P^n - \pi\|_{\text{TV}}$ and $\|\xi P^n - \pi\|_f$, respectively.

To compare our result with Meyn and Tweedie (1994), Rosenthal (1995) and Roberts and Tweedie (1999), we now derive from the explicit expressions of the bounds provided in Theorem 2 the rate of convergence for the total variation distance or the $f$-norm, that is, we find a bound for $\limsup_{n \to \infty} n^{-1} \log \|P^n(x, \cdot) - \pi\|_f$. We follow the approach originally taken by Rosenthal (1995), but we adapt the results to the expression of the bound given in Theorem 2.

PROPOSITION 3. *Assume* (A1) *and* (A2), *and that* $\pi P = \pi$. *Let* $f : \mathcal{X} \to [1, \infty)$ *be a function that satisfies* $f(x) + f(x') \leq 2\bar{V}(x, x')$ *for all* $(x, x') \in \mathcal{X} \times \mathcal{X}$. *Then, for all* $x \in \mathcal{X}$,

$$
\limsup_{n \to \infty} n^{-1} \log \|P^n(x, \cdot) - \pi\|_f
$$
(18)
$$
\leq \begin{cases} \dfrac{-\log(\lambda) \log(1 - \varepsilon)}{\log((M - \varepsilon)/\lambda) - \log(1 - \varepsilon)}, & \text{if } \dfrac{M - \varepsilon}{\lambda} \geq 1, \\ \log(\lambda), & \text{if } \dfrac{M - \varepsilon}{\lambda} < 1, \end{cases}
$$

*where* $M := \sup_{(x, x') \in \bar{C}} \bar{P}\bar{V}(x, x')$.

PROOF. By definition of $\bar{P}$ [see (3)], for all $(x, x') \in \bar{C}$ we have

$$(1 - \varepsilon)\bar{R}\bar{V}(x, x') + \varepsilon \int \nu_{x,x'}(dy) \bar{V}(y, y) = \bar{P}\bar{V}(x, x') \geq (1 - \varepsilon)\bar{R}\bar{V}(x, x') + \varepsilon,$$

where we have used that $\bar{V} \geq 1$. Thus

$$\sup_{(x,x') \in \bar{C}} \bar{R}\bar{V}(x, x') \leq \frac{M - \varepsilon}{1 - \varepsilon},$$

which implies

$$(1 - \varepsilon) \sup_{(x,x') \in \bar{C}} \bar{R}\bar{V}(x, x') \lambda^{-1} \leq (M - \varepsilon)\lambda^{-1}.$$

Assuming first $(M - \varepsilon)\lambda^{-1} \geq 1$, the bounds for total variation and $f$-norm can be expressed for $j \in \{1, \ldots, n\}$,

$$\|P^n(x, \cdot) - \pi\|_{\text{TV}} \leq 2(1 - \varepsilon)^j + 2\lambda^{n-j+1}(M - \varepsilon)^{j-1} \int \bar{V}(x, x') \pi(dx'),$$

$$\|P^n(x, \cdot) - \pi\|_f \leq \frac{2b(1 - \varepsilon)^j}{1 - \lambda}$$
$$+ 2\lambda^n ((1 - \varepsilon)^j + \lambda^{-j+1}(M - \varepsilon)^{j-1}) \int \bar{V}(x, x') \pi(dx').$$

The result follows by choosing
$$j = \left\lfloor \frac{-\log(\lambda)n}{\log((M-\varepsilon)/\lambda) - \log(1-\varepsilon)} \right\rfloor.$$

When $(M-\varepsilon)\lambda^{-1} < 1$, we put $j = n+1$ in (11) and (12), showing that
$$\|\xi P^n - \xi' P^n\|_{\text{TV}} \leq 2\lambda^n \int \bar{V}(x,x')\pi(dx') \quad \text{and}$$
$$\|\xi P^n - \xi' P^n\|_f \leq 2\lambda^n \int \bar{V}(x,x')\pi(dx').$$

The result follows. □

REMARK 1. The bounds we find in this paper for the $f$-total variation distance are the same as those found for the total variation distance by Roberts and Tweedie [(1999), Theorem 2.3].

In some applications, the minorization and drift conditions (A1) and (A2) are more naturally expressed in terms of the kernel $P$, and it is thus required to derive the bivariate drift and minorization conditions from the corresponding single variate conditions [Rosenthal (1995), Theorem 12, and Roberts and Tweedie (1999), Section 5]. The crucial point here is to relate the bivariate drift condition (A2) to single variate drift condition. We essentially follow to Rosenthal's [(1995), Theorem 12] argument, which allows us to construct such a drift function $\bar{V}$ from univariate test functions [see Roberts and Tweedie (1999), Theorem 5.2, for a refinement of this result].

Consider the following assumption:

(S) There exist a function $V$ and a constant $c$ such that:

- The level set $C = \{x \in \mathcal{X} : V(x) \leq c\}$ is $(1,\varepsilon)$-small; that is, $P(x,\cdot) \geq \varepsilon\nu(\cdot)$ for all $x \in C$ for some $\varepsilon > 0$ and some probability measure $\nu$.
- There exist $\lambda_c < 1$ and $b_c < \infty$ such that $PV \leq \lambda_c V + b_c \mathbf{1}_C$ and $\lambda_c + b_c/(1+c) < 1$.

Under (S), $\bar{C} = \{(x,x'); V(x) \leq c, V(x') \leq c\}$ is a $(1,\varepsilon)$-coupling set, that is, for all $(x,x') \in \bar{C}$ and all $A \in \mathcal{B}(\mathcal{X})$, $P(x,A) \wedge P(x',A) \geq \varepsilon\nu(A)$. Define the univariate residual kernel $R$ as

(19) $\quad R(x,A) = (1-\varepsilon)^{-1}(P(x,A) - \varepsilon\nu(A)) \qquad \forall x \in C, \ \forall A \in \mathcal{B}(\mathcal{X}).$

To apply Theorem 1, we need to define the kernels $\bar{R}$, $\bar{P}$ and $P^*$. Because the drift condition is expressed on the univariate kernel $P$, we define both $\bar{R}$ and $\bar{P}$ from the corresponding univariate kernels $R$ and $P$. More precisely, for all $A, A' \in \mathcal{B}(\mathcal{X})$, define

(20) $\quad \bar{R}(x,x'; A \times A') := R(x,A)R(x',A') \qquad \text{if } (x,x') \in \bar{C},$

(21) $\quad \bar{P}(x,x'; A \times A') := P(x,A)P(x',A') \qquad \text{if } (x,x') \notin \bar{C}.$

These kernels satisfy (2) and (4).

PROPOSITION 4. *Assume* (S). *Then* (A1) *is satisfied with* $\bar{C} = C \times C$ *and* $\nu_{x,x'} = \nu$ *for all* $(x, x') \in C \times C$. *Define* $P^*$ *as in* (6) *with* $\bar{R}$ *and* $\bar{P}$ *given in* (20) *and* (21). *Then* (A2) *is satisfied with* $\bar{V}(x, x') = (1/2)(V(x) + V(x'))$ *for all* $(x, x') \in \mathcal{X} \times \mathcal{X}$ *with*

$$\lambda = \lambda_c + b_c/(1+c) \quad and \quad b = \left\{\frac{c\varepsilon\lambda_c}{1-\varepsilon} - \frac{cb_c}{1+c}\right\} \vee 0 + \frac{b_c - \varepsilon}{1-\varepsilon}.$$

PROOF. The proof follows from Roberts and Tweedie [(1999), Theorem 5.2]. Since, for $(x, x') \notin \bar{C}$, $(1+c)/2 \leq \bar{V}(x, x')$, we have

$$P^*\bar{V}(x, x') \leq \lambda_c \bar{V}(x, x') + \frac{b_c}{2} \leq \left(\lambda_c + \frac{b_c}{1+c}\right)\bar{V}(x, x') \qquad \forall\, (x, x') \notin C \times C$$

and, for $(x, x') \in C \times C$,

$$\begin{aligned}
P^*\bar{V}(x, x') &= \frac{1}{2}(RV(x) + RV(x')) \\
&= \frac{1}{2(1-\varepsilon)}(PV(x) + PV(x') - 2\varepsilon\nu(V)) \\
&\leq \frac{\lambda_c}{(1-\varepsilon)}\bar{V}(x, x') + \frac{b_c - \varepsilon}{1-\varepsilon} \leq \lambda \bar{V}(x, x') + b,
\end{aligned}$$

where we have used that, for $(x, x') \in \bar{C}$, $\bar{V}(x, x') \leq c$. The proof follows. □

Under (S), we may thus apply Theorem 2 with $f = V$ which yields explicit bounds for the total variation and the $V$-norm, under the assumptions used by Rosenthal (1995) and Roberts and Tweedie (1999) to obtain bounds for the total variation distance [see also Rosenthal (2002)]. It is worthwhile to note that (see the discussion above) the rate of convergence in $V$-norm is the *same* as the rate of convergence in total variation.

REMARK 2. It may be checked that if the sets $\{V \leq d\}$ are 1-small for all $d \geq c$, then assumption (S) is always satisfied for large enough $d$ [see Roberts and Tweedie (1999), discussion following Theorem 5.2].

We summarize the discussion above in the following theorem.

THEOREM 5. *Assume* (S). *Then, for all* $j \in \{1, \ldots, n+1\}$ *and for all initial probability measures* $\xi$ *and* $\xi'$ *on* $\mathcal{X}$,

$$\|\xi P^n - \xi' P^n\|_{\mathrm{TV}} \leq 2(1-\varepsilon)^j \mathbf{1}(j \leq n) + \lambda^n B^{j-1}(\xi(V) + \xi'(V)),$$

$$\begin{aligned}
\|\xi P^n - \xi' P^n\|_V &\leq 2(1-\varepsilon)^j (b(1-\lambda)^{-1} + \lambda^n(\xi(V) + \xi'(V))/2)\mathbf{1}(j \leq n) \\
&\quad + \lambda^n B^{j-1}(\xi(V) + \xi'(V)),
\end{aligned}$$

where $\lambda = \lambda_c + b_c/(1+c)$ and

$$B = 1 \vee \left((1-\varepsilon)\lambda^{-1} \sup_{x \in C} RV(x)\right).$$

1.5. *Example.* We conclude this section with a simple example that shows a situation where we can exploit the additional degree of flexibility brought by $(1,\varepsilon)$-coupling sets. Consider the Markov chain on $\mathbb{R}^d$ defined for $k \in \mathbb{Z}^+$ by

$$X_{k+1} = g(X_k) + Z_k,$$

where:

1. $g$ is a Lipshitz function over $\mathbb{R}^d$ for some norm $\|\cdot\|$ with Lipshitz constant

$$\|g\|_{\mathrm{Lip}} = \sup_{\substack{(x,y)\in\mathbb{R}^d\times\mathbb{R}^d \\ x\neq y}} \frac{\|g(x)-g(y)\|}{\|x-y\|} < 1;$$

2. $(Z_k, k \geq 0)$ is a sequence of independent and identically distributed random vectors with density $q$ w.r.t. Lebesgue measure on $\mathbb{R}^d$. In addition, $q$ is positive and continuous.

It is known [see, e.g., Doukhan and Ghindes (1980)] that under these assumptions the Markov chain is positive recurrent and thus has a unique invariant distribution. Define for $\delta > 0$,

(22) $$\bar{C}(\delta) := \{(x,x') \in \mathbb{R}^d \times \mathbb{R}^d : \|x - x'\| \leq \delta\}.$$

Using $a \wedge b = (1/2)((a+b) - |a-b|)$, it is easily shown that for all $(x,x') \in \bar{C}(\delta)$ and all $A \in \mathcal{B}(\mathbb{R}^d)$,

$$P(x,A) \wedge P(x',A)$$
$$\geq \tfrac{1}{2} \int_A (q(z-g(x)) + q(z-g(x'))) - |q(z-g(x)) - q(z-g(x'))|) \, dz$$

and thus $P(x,A) \wedge P(x',A) \geq \varepsilon(\delta)\nu_{x,x'}(A)$ with

$$\nu_{x,x'}(A) = \int_A \Big(q(z-g(x)) + q(z-g(x')) $$
$$- |q(z-g(x)) - q(z-g(x'))|\Big) dz$$

(23)
$$\times \left(2 - \int |q(z-g(x)) - q(z-g(x'))| \, dz\right)^{-1},$$

$$\varepsilon(\delta) = 1 - \tfrac{1}{2} \sup_{(x,x')\in\bar{C}(\delta)} \int |q(z-(g(x)-g(x'))) - q(z)| \, dz.$$

Note that for all $(x,x') \in \bar{C}(\delta)$, $\|g(x) - g(x')\| \leq \|g\|_{\text{Lip}} \|x - x'\| \leq \|g\|_{\text{Lip}} \delta$. Since the function $u \to \int |q(z-u) - q(z)| \, dz$ is continuous and $q$ is everywhere positive, for all $\delta > 0$, the set $\bar{C}(\delta)$ is a $(1, \varepsilon(\delta))$-coupling set.

Let $\delta > 0$. For all $(x,x') \in \mathbb{R}^d \times \mathbb{R}^d$ and all $A, A' \in \mathcal{B}(\mathbb{R}^d)$, define $\bar{P}$ by

$$\bar{P}(x,x'; A \times A') = \int \mathbf{1}_A(f(x) + z) \mathbf{1}_{A'}(f(x') + z) q(z) \, dz$$

and let, for $(x,x') \in \bar{C}(\delta)$,

$$\bar{R}_\delta(x,x'; A \times A') = (1 - \varepsilon(\delta))^{-1} (\bar{P}(x,x'; A \times A') - \varepsilon(\delta) \nu_{x,x'}(A \cap A')).$$

It is easily checked that $\bar{R}_\delta$ and $\bar{P}$ satisfy (2) and (4), respectively. Finally, define $P_\delta^*$ as in (6).

We now determine an explicit bound for the total variation distance. Put $\bar{V}(x,x') = 1 + \|x - x'\|$. Note that for all $(x,x') \in \mathbb{R}^d \times \mathbb{R}^d$,

$$\bar{P}\bar{V}(x,x') = 1 + \|g(x) - g(x')\| \leq 1 + \|g\|_{\text{Lip}} \|x - x'\|.$$

Choose $\lambda$ such that $\|g\|_{\text{Lip}} < \lambda < 1$. By construction, for all $(x,x') \notin \bar{C}(\delta)$, we have $\|x - x'\| \geq \delta$. Hence, for any $\delta > (1-\lambda)/(\lambda - \|g\|_{\text{Lip}})$ and all $(x,x') \notin \bar{C}(0,\delta)$, we have

$$1 + \|g\|_{\text{Lip}} \|x - x'\| = \lambda(1 + \|x - x'\|) + (1 - \lambda - (\lambda - \|g\|_{\text{Lip}}) \|x - x'\|)$$
$$\leq \lambda(1 + \|x - x'\|) + (1 - \lambda - (\lambda - \|g\|_{\text{Lip}})\delta)$$
$$< \lambda(1 + \|x - x'\|).$$

It remains to prove that $\sup_{(x,x') \in \bar{C}(\delta)} \bar{R}\bar{V}(x,x') < \infty$. Note that

$$\sup_{(x,x') \in \bar{C}(\delta)} \bar{R}\bar{V}(x,x') \leq \frac{\sup_{(x,x') \in \bar{C}(\delta)} \bar{P}\bar{V}(x,x') - \varepsilon(\delta)}{1 - \varepsilon(\delta)} \leq \frac{1 + \|g\|_{\text{Lip}} \delta - \varepsilon(\delta)}{1 - \varepsilon(\delta)}.$$

Summarizing our findings, for any $\lambda$ with $\|g\|_{\text{Lip}} < \lambda < 1$ and any $\delta > (1-\lambda)/(\lambda - \|g\|_{\text{Lip}})$, (A1) is satisfied with $\varepsilon := \varepsilon(\delta)$ and (A2) is satisfied with $\bar{V}(x,x') = 1 + \|x - x'\|$. We may thus apply Theorem 2 to obtain a total variation distance bound as follows. [Note that with this choice of bivariate drift function $\bar{V}$ we may only compute total variation bound; the condition $f(x) + f(x') \leq 2(1 + \|x - x'\|)$ indeed implies that $f \leq 1$.]

PROPOSITION 6. *For all $\lambda$ such that $\|g\|_{\text{Lip}} < \lambda < 1$, for all $\delta > (1-\lambda)/(\lambda - \|g\|_{\text{Lip}})$, for all $j \in \{1, \ldots, n+1\}$ and for all initial probability measures $\xi$ and $\xi'$ on $\mathcal{X}$,*

$$\|\xi P^n - \xi' P^n\|_{\text{TV}}$$
$$\leq 2(1 - \varepsilon(\delta))^j \mathbf{1}(j \leq n) + 2\lambda^n B^{j-1} \left(1 + \int\int \xi(dx) \xi'(dx') \|x - x'\|\right),$$

*where $\varepsilon(\delta)$ is defined in (23) and*

$$B = 1 \vee \{\lambda^{-1}(1 + \|g\|_{\text{Lip}} \delta - \varepsilon(\delta))\}.$$

**2. Time-inhomogeneous case.** We now proceed to extend Theorem 2 to time-inhomogeneous chains. Specifically, we consider a family $(P_k, k \geq 1)$ of Markov transition kernels. That is, we allow $P_k(x, A)$ to depend not only on the starting point $x$ and the target subset $A$, but also on the time parameter $k$. For example, this would be the case for simulated annealing and hidden Markov models; a specific example is discussed in Section 3.

2.1. *Assumptions and lemma.* The assumptions and notations parallel those from the time-homogeneous case. We first assume the following minorization condition.

(NS1) There exist a sequence $(\bar{C}_k, k \geq 1)$ of subsets of $\mathcal{X} \times \mathcal{X}$, $\bar{C}_k \subset \mathcal{X} \times \mathcal{X}$, a sequence $(\varepsilon_k, k \geq 1)$, $\varepsilon_k \geq 0$, and a family of probability measures $(\nu_{k,x,x'}, (x, x') \in \bar{C}_k, k \geq 1)$ such that

$$P_k(x, \cdot) \wedge P_k(x', \cdot) \geq \varepsilon_k \nu_{k,x,x'}(\cdot).$$

Let $(\bar{P}_k, k \geq 1)$ be a family of transitions kernels that satisfy, for all $k$, the analog of (4) with $P = P_k$ and let $(\bar{R}_k, k \geq 0)$ be a family of transition kernels that verify, for all $k$, the analog of (2) with $P = P_k$, $\nu_{x,x'} = \nu_{k,x,x'}$, $\varepsilon = \varepsilon_k$ and $\bar{C} = \bar{C}_k$. The proof is based on straightforward adaptation of the coupling construction used in the homogeneous case. For $n \geq 0$, if $(X_n, X_n') \in \bar{C}_{n+1}$ and $d_n = 0$, flip a coin with probability of success $\varepsilon_{n+1}$. If the coin comes up heads, then draw $X_{n+1}$ from $\nu_{n+1, X_n, X_n'}$ and set $X_{n+1} = X_{n+1}'$ and $d_{n+1} = 1$. Otherwise, draw $(X_{n+1}, X_{n+1}')$ from $\bar{R}_{n+1}(X_n, X_n'; \cdot)$ and set $d_{n+1} = 0$. If $(X_n, X_n') \notin \bar{C}_{n+1}$ and $d_n = 0$, then draw $(X_{n+1}, X_{n+1}')$ from $\bar{P}_{n+1}(X_n, X_n'; \cdot)$ and set $d_{n+1} = 0$. Finally, define $(P_k^*, k \geq 0)$ to be the family of transition kernels defined as the analog of (6). For $\mu$ a probability measure on $\mathcal{X} \times \mathcal{X}$, denote $\mathbb{P}_\mu^*$ and $\mathbb{E}_\mu^*$ the probability and the expectation induced by the Markov chain with initial distribution $\mu$ and transition kernels $(P_k^*, k \geq 0)$.

LEMMA 7. *Assume* (NS1) *and let* $f : \mathcal{X} \to [1, +\infty)$. *For any probability measures* $\xi$, $\xi'$ *on* $\mathcal{B}(\mathcal{X})$, *for any* $n \geq 1$,

$$
\begin{aligned}
&\|\xi P_1 \cdots P_n - \xi' P_1 \cdots P_n\|_f \\
&\qquad \leq \mathbb{E}_{\xi \otimes \xi'}^* \bigg\{ (f(X_n) + f(X_n')) \prod_{i=1}^n (1 - \varepsilon_i \mathbf{1}_{\bar{C}_i}(\bar{X}_{i-1})) \bigg\},
\end{aligned}
\tag{24}
$$

*where* $\bar{X}_i = (X_i, X_i')$.

The proof can be adapted from Lemma 1 and (9). We also assume the following drift condition:

(NS2) There exist a family of functions $\{\bar{V}_k\}_{k\geq 0}$, $\bar{V}_k : \mathcal{X} \times \mathcal{X} \to [1,\infty)$, and two sequences $(\lambda_k, k \geq 0)$, $0 \leq \lambda_k \leq 1$ for all $k \geq 1$ and $(b_k, k \geq 0)$, such that

$$P_{k+1}^* \bar{V}_{k+1} \leq \lambda_k \bar{V}_k + b_k \mathbf{1}_{\bar{C}_{k+1}} \qquad \forall k \geq 0. \tag{25}$$

Define for $j \in \{1, \ldots, k\}$,

$$(1-\varepsilon)_{j,k} := \max_{1 \leq k_1 < \cdots < k_j \leq k} \prod_{l=1}^{j} (1-\varepsilon_{k_l}) \quad \text{and} \quad B_{j,k} := \max_{1 \leq k_1 < \cdots < k_j \leq k} \prod_{l=1}^{j} B_{k_l},$$

where, for any integer $k$,

$$B_k := 1 \vee \left( (1-\varepsilon_k) \left( \sup_{(x,x') \in \bar{C}_k} \bar{R}_k \bar{V}_k(x,x') \right) \lambda_{k-1}^{-1} \right). \tag{26}$$

By convention, we set $B_{0,k} = 1$ for all $k$.

2.2. *Main time-inhomogeneous result.* We can now state our main result, as follows.

THEOREM 8. *Assume* (NS1) *and* (NS2). *Let* $(f_k, k \geq 0)$ *be a family of functions such that, for all* $k \geq 0$, $f_k(x) + f_k(x') \leq 2\bar{V}_k(x,x')$ *for all* $(x,x') \in \mathcal{X} \times \mathcal{X}$. *Then, for all* $j \in \{1, \ldots, n+1\}$ *and for all initial probability measures* $\xi$ *and* $\xi'$,

$$\|\xi P_1 \cdots P_n - \xi' P_1 \cdots P_n\|_{\mathrm{TV}} \tag{27}$$
$$\leq 2(1-\varepsilon)_{j,n} \mathbf{1}(j \leq n) + 2 \left( \prod_{s=0}^{n-1} \lambda_s \right) B_{j-1,n} (\xi \otimes \xi')(\bar{V}_0),$$

$$\|\xi P_1 \cdots P_n - \xi' P_1 \cdots P_n\|_{f_n} \tag{28}$$
$$\leq 2(1-\varepsilon)_{j,n} D_n \mathbf{1}(j \leq n) + 2 \left( \prod_{s=0}^{n-1} \lambda_s \right) B_{j-1,n} (\xi \otimes \xi')(\bar{V}_0),$$

where $D_n := (\prod_{l=0}^{n-1} \lambda_l) \xi \otimes \xi'(V_0) + \sum_{j=0}^{n-1} (\prod_{l=j+1}^{n-1} \lambda_l) b_j$ with the convention $\prod_{l=i}^{j} \lambda_l = 1$ when $i > j$.

PROOF. The proof is along the same lines as for the time-homogeneous case. Denote $N_k = \sum_{j=0}^{k} \mathbf{1}_{\bar{C}_{j+1}}(X_j, X_j')$. For any $j \in \{1, \ldots, n+1\}$, we have

$$\mathbb{E}_{\xi \otimes \xi'}^* \left\{ (f_n(X_n) + f_n(X_n')) \prod_{i=1}^{n} (1 - \varepsilon_i \mathbf{1}_{\bar{C}_i}(\bar{X}_{i-1})) \right\}$$
$$\leq (1-\varepsilon)_{j,n} \mathbb{E}_{\xi \otimes \xi'}^* \{ (f_n(X_n) + f_n(X_n')) \} \mathbf{1}(j \leq n)$$
$$+ 2 \mathbb{E}_{\xi \otimes \xi'}^* \left\{ \bar{V}_n(\bar{X}_n) \prod_{i=1}^{n} (1 - \varepsilon_i \mathbf{1}_{\bar{C}_i}(\bar{X}_{i-1})) \mathbf{1}(N_{n-1} < j) \right\},$$

where we have used that $\prod_{i=1}^{n}(1 - \varepsilon_i \mathbf{1}_{\bar{C}_i}(\bar{X}_{i-1}))\mathbf{1}(N_{n-1} \geq j) \leq (1 - \varepsilon)_{j,n}$. When $f_n \equiv 1$,

$$\mathbb{E}^*_{\xi \otimes \xi'}\{(f_n(X_n) + f_n(X'_n))\} = 2.$$

Otherwise,

$$\mathbb{E}^*_{\xi \otimes \xi'}\{(f_n(X_n) + f_n(X'_n))\} \leq 2\mathbb{E}^*_{\xi \otimes \xi'}\{\bar{V}_n(\bar{X}_n)\} \leq 2D_n.$$

Now, since by definition $B_j \geq 1$ [see (26)], we have $B_{j,n} \leq B_{j',n}$ for all $0 \leq j \leq j' \leq n$ and

$$\mathbf{1}(N_{n-1} \leq j - 1)(B_{j-1,n})^{-1} \leq (B_{N_{n-1},n})^{-1},$$

which implies that

$$(29) \quad \mathbb{E}^*_{\xi \otimes \xi'}\left\{\bar{V}_n(\bar{X}_n) \prod_{i=1}^{n}(1 - \varepsilon_i \mathbf{1}_{\bar{C}_i}(\bar{X}_{i-1}))\mathbf{1}(N_{n-1} < j)\right\}$$
$$\leq \left(\prod_{j=0}^{n-1} \lambda_j\right) B_{j-1,n} \mathbb{E}^*_{\xi \otimes \xi'}\{M_n\},$$

where, for $s \geq 0$,

$$(30) \quad M_s := \left(\prod_{j=0}^{s-1} \lambda_j\right)^{-1}(B_{1,N_{s-1},s})^{-1} \prod_{j=1}^{s}(1 - \varepsilon_j \mathbf{1}_{\bar{C}_j}(\bar{X}_{j-1}))\bar{V}_s(X_s, X'_s).$$

As above, $(M_s, s \geq 0)$ is an $(\mathcal{F}, \mathbb{P}^*_{\xi \otimes \xi'})$ supermartingale w.r.t., where $\mathcal{F} := \{\mathcal{F}_s := \sigma(\bar{X}_j, 0 \leq j \leq s), s \geq 0\}$, which concludes the proof. □

**3. Application to simulated annealing.** In this section, we apply the results above to study the convergence of the simulated annealing (SA) algorithm for continuous global optimization [see Locatelli (2001, 2002), Fouskakis and Draper (2001), Andrieu, Breyer and Doucet (2001) and the references therein].

3.1. *Assumptions.* Let $f$ be a function defined on $\mathbb{R}$, and let $\mathcal{M}$ be the set of global minima of $f$ (to keep the discussion simple, multidimensional versions are not considered here). We make the following assumptions:

(SA0) the function $f$ is twice continuously differentiable and there exist $\alpha > 0$, $x_1 \in \mathbb{R}$, such that, for all $y \geq x \geq x_1$,

(31) $$f(y) - f(x) \geq \alpha(y - x)$$

and similarly, for all $y \leq x \leq -x_1$,

(32) $$f(y) - f(x) \geq \alpha(x - y).$$

(SA1) For each $x \in \mathcal{M}$, we have $f''(x) > 0$.

Under (SA0), $\mathcal{M} \subseteq [-x_1, x_1]$, that is, the set of global minima of $f$ is contained in the interval $[-x_1, x_1]$. Assumption (SA1) implies that the global minima are isolated and thus, that the set $\mathcal{M}$ is finite. Assumption (SA0) implies that for all $\gamma \geq 0$, $\int \exp(-\gamma f(y)) \mu^{\text{Leb}}(dy) < \infty$, where $\mu^{\text{Leb}}$ is the Lebesgue measure over $\mathbb{R}$.

Consider a *candidate transition kernel*, $Q(x, A)$, $x \in \mathbb{R}$, $A \in \mathcal{B}(\mathbb{R})$, which generates potential transitions for a discrete time Markov chain evolving on $\mathbb{R}$. We focus on the case where the candidate points are proposed from a random walk with increment distribution that has a density $q$ with respect to $\mu^{\text{Leb}}$: $Q(x, A) = \int_A q(y - x) \mu^{\text{Leb}}(dy)$, $A \in \mathcal{B}(\mathbb{R})$. In addition, make the following assumption:

(SA2) The proposal density $q$ is continuous and strictly positive and symmetric: $q(y) > 0$ and $q(y) = q(-y)$.

3.2. *The random walk Metropolis–Hastings algorithm.* The random walk Metropolis–Hastings (RWMH) algorithm corresponds to the Hastings–Metropolis algorithm introduced by Metropolis, Rosenbluth, Rosenbluth, Teller and Teller (1953) and Hastings (1970). It proceeds as follows to sample from the (unnormalized) distribution $\exp(-\gamma f(x)) \mu^{\text{Leb}}(dx)$ for $\gamma > 0$. (For RWMH, the "inverse temperature" parameter $\gamma$ is held constant. We see later that with simulated annealing, by contrast, $\gamma$ is modified at each iteration of the algorithm.)

Given the current state $x$, a candidate new state $y$ is chosen according to the law $Q(x, \cdot)$. This candidate $y$ is then accepted with probability $\alpha_\gamma(x, y)$, where

$$\alpha_\gamma(x, y) = 1 \wedge (\exp(-\gamma(f(y) - f(x)))).$$

The RWMH kernel is thus given by

$$
\begin{aligned}
K_\gamma(x, A) &= \int_A \alpha_\gamma(x, y) q(y - x) \mu^{\text{Leb}}(dy) \\
&\quad + \delta_x(A) \int (1 - \alpha_\gamma(x, y)) q(y - x) \mu^{\text{Leb}}(dy), \qquad A \in \mathcal{B}(\mathbb{R}).
\end{aligned}
$$
(33)

It then follows that $\pi_\gamma(\cdot)$ is a stationary distribution for $K_\gamma$, where

$$\pi_\gamma(A) = \frac{\int_A \exp(-\gamma f(x)) \mu^{\text{Leb}}(dx)}{\int_{\mathbb{R}} \exp(-\gamma f(x)) \mu^{\text{Leb}}(dx)} \qquad \forall A \in \mathcal{B}(\mathbb{R}).$$

The RWMH algorithm on $\mathbb{R}$ was extensively studied by Mengersen and Tweedie (1996), who showed that the transition kernels $K_\gamma$ are $\pi_\gamma$-irreducible (Lemma 1.1) and that all the compact sets are small (Lemma 1.2).

LEMMA 9. *Assume* (SA0)–(SA2). *Then, for every compact subset $C$ of $\mathbb{R}$ such that $\mu^{\mathrm{Leb}}(C) > 0$, we have for all $x \in C$, $K_\gamma(x, A) \geq \varepsilon_\gamma \nu_\gamma(A)$ with*

$$(34) \qquad \varepsilon_\gamma := \varepsilon e^{-\gamma d} \lambda^{\mathrm{Leb}}(C) \quad and \quad \nu(A) := \frac{\lambda^{\mathrm{Leb}}(A \cap C)}{\lambda^{\mathrm{Leb}}(C)},$$

*where*

$$(35) \qquad d := \sup_{x \in C} f(x) - \inf_{x \in C} f(x) \quad and \quad \varepsilon := \inf_{(x,y) \in C \times C} q(y - x) > 0.$$

PROOF. For all $x \in C$,

$$K_\gamma(x, A) \geq \int_{A \cap C} (e^{-\gamma(f(y) - f(x))} \wedge 1) q(y - x) \mu^{\mathrm{Leb}}(dy) \geq \varepsilon e^{-\gamma d} \lambda^{\mathrm{Leb}}(A \cap C).$$

□

To apply Theorem 8, we need to find drift functions that satisfy drift conditions outside the compact sets of $\mathbb{R}$. The existence of drift functions for the RWMH algorithm was shown by Mengersen and Tweedie [(1996), Theorem 3.2]. The proposition below relaxes some of the assumptions required in their result, and shows that the same drift function can be taken for all the Markov kernels $K_\gamma$ for large enough $\gamma$. For $0 < s \leq \gamma$, let $V_s(x) := e^{sf(x)}$ and

$$(36) \qquad r(\gamma, s) := 1 - \left(\frac{\gamma - s}{\gamma}\right)^{\gamma/s} + \left(\frac{\gamma - s}{\gamma}\right)^{(\gamma - s)/s}.$$

PROPOSITION 10. *Assume* (SA0)–(SA2). *Then, for all $\beta$ such that $1/2 < \beta < 1$, there exist $\underline{x} < \infty$, $\underline{\gamma} > 0$ and $s > 0$ such that:*

(i) $(K_\gamma V_s(x))/(V_s(x)) \leq r(\gamma, s)$ *for all $x \in \mathbb{R}$ and $\gamma \geq 0$;*
(ii) $(K_\gamma V_s(x))/(V_s(x)) \leq \beta$ *for all $|x| \geq \underline{x}$ and $\gamma \geq \underline{\gamma}$.*

PROOF. By (33) and using that $V_s(y) = e^{sf(x)}$, we have, for $\gamma > s > 0$,

$$(37) \qquad \frac{K_\gamma V_s(x)}{V_s(x)} = \int \varphi_{\gamma,s}(e^{-(f(y) - f(x))}) q(y - x) \mu^{\mathrm{Leb}}(dy),$$

where $\varphi_{\gamma,s}(u) := u^{-s}(u^\gamma \wedge 1) + 1 - (u^\gamma \wedge 1)$. We easily check that, for all $u \geq 0$,

$$(38) \qquad \varphi_{\gamma,s}(u) \leq \varphi_{\gamma,s}\left[\left(\frac{\gamma - s}{\gamma}\right)^{1/s}\right] = r(\gamma, s),$$

which proves the first assertion of the proposition. Now, for any $\varepsilon > 0$, we prove that there exists some $\underline{x}$ such that

$$\lim_{\gamma \to \infty} \sup_{x \geq \underline{x}} \frac{K_\gamma V_s(x)}{V_s(x)} \leq \varepsilon + \frac{1}{2}.$$

The proof of the corresponding inequality where $x \geq \underline{x}$ is replaced by $x \leq -\underline{x}$ follows the same lines. Choose $M > 0$ such that

$$\int_{-\infty}^{-M} q(z) \mu^{\mathrm{Leb}}(dz) \leq \varepsilon/2.$$

Inserting this inequality into (37), where $z = y - x$, and using (38) yields

$$\frac{K_\gamma V_s(x)}{V_s(x)} \leq \int_{-M}^{0} \varphi_{\gamma,s}(e^{-(f(x+z)-f(x))}) q(z) \mu^{\mathrm{Leb}}(dz)$$
$$+ r(\gamma, s) \left( \int_0^\infty q(z) \mu^{\mathrm{Leb}}(dz) + \frac{\varepsilon}{2} \right).$$

For all $x \geq \underline{x} := x_1 + M$ and all $-M \leq z \leq 0$, we have by assumption (SA0), $\exp(-(f(x+z) - f(x))) \geq \exp(-\alpha z) \geq 1$ and since $\varphi_{\gamma,s}(u) = u^{-s}$ for $u \geq 1$,

$$\frac{K_\gamma V_s(x)}{V_s(x)} \leq \int_{-M}^{0} e^{\alpha s z} q(z) \mu^{\mathrm{Leb}}(dz) + r(\gamma, s) \frac{1+\varepsilon}{2}.$$

Now, choose $s$ sufficiently large so that the first term on the right-hand side is less than $\varepsilon/2$. Once $s$ is chosen, we easily check that $\lim_{\gamma \to \infty} r(\gamma, s) = 1$. This proves the second assertion. □

Define $\bar{K}_\gamma(x, x'; A \times A') = K_\gamma(x, A) K_\gamma(x', A')$ and, for $s \geq 0$, $\bar{V}_s(x, x') = (1/2)(V_s(x) + V_s(x'))$.

PROPOSITION 11. *Assume* (SA0)–(SA2). *For all $s \geq 0$ and for all $c \geq 0$, $\{V_s \leq c\}$ is a compact 1-small set for $K_\gamma$. Moreover, there exist $0 \leq \lambda_0 < \lambda < 1$, $s > 0$, $c_0 \leq c$, $b$ and $\underline{\gamma}$ such that, for all $\gamma \geq \underline{\gamma}$,*

(39) $$K_\gamma V_s \leq \lambda_0 V_s + b \mathbf{1}_{\{V_s \leq c_0\}},$$

(40) $$\bar{K}_\gamma \bar{V}_s \leq \lambda \bar{V}_s + b \mathbf{1}_{\{V_s \leq c\} \times \{V_s \leq c\}}.$$

PROOF. The compactness of $\{V_s \leq c\}$ is straightforward from (SA0). Then, by Lemma 9, it is a 1-small set for $K_\gamma$. Equation (39) follows from Proposition 10. To prove (40), write for $c \geq c_0$,

$$\bar{K}_\gamma \bar{V}_s \leq \lambda_0 \bar{V}_s + b \mathbf{1}_{\{V_s \leq c\} \times \{V_s \leq c\}} + (b/2)(\mathbf{1}_{\{V_s \leq c\} \times \{V_s > c\}} + \mathbf{1}_{\{V_s > c\} \times \{V_s \leq c\}}).$$

Set $0 \leq \lambda_0 < \lambda < 1$ and $c = (b/(\lambda - \lambda_0) - 1) \vee c_0$. We have, for all $(x, x') \in \{V_s \leq c\} \times \{V_s > c\}$,

$$b/2 \leq (\lambda - \lambda_0)(1+c)/2 \leq (\lambda - \lambda_0) \bar{V}_s(x, x'),$$

which implies

$$(\lambda_0 \bar{V}_s + (b/2))\mathbf{1}_{\{V_s \leq c\} \times \{V_s > c\}} \leq \lambda \bar{V}_s \mathbf{1}_{\{V_s \leq c\} \times \{V_s > c\}}.$$

This concludes the proof. □

The key point in the above result [also outlined in Andrieu, Breyer and Doucet (2001)] is that, for large enough $\gamma$ ($\gamma \geq \underline{\gamma}$), all the transition kernels $\bar{K}_\gamma$ satisfy a drift condition outside the *same small set* $\{V_s \leq c\} \times \{V_s \leq c\}$, with the *same* drift function $\bar{V}_s$ and the *same* constants $\lambda$ and $b$.

3.3. *The simulated annealing algorithm.* We now consider the simulated annealing case. Here $\gamma = \gamma_i$ depends on the iteration, and for the $i$th iteration, the kernel $P_i = K_{\gamma_i}$ is used. Define similarly $\bar{P}_i = \bar{K}_{\gamma_i}$ and $\pi_i = \pi_{\gamma_i}$. Denote $\bar{C} = \{V_s \leq c\} \times \{V_s \leq c\}$, with the constants $s$ and $c$ chosen to satisfy (40). For $(x, x') \in \bar{C}$, set $\bar{R}_i(x, x'; A \times A') = R_i(x, A)R_i(x', A')$, with

(41) $\quad R_i(x, A) = (1 - \varepsilon_i)^{-1}(P_i(x, A) - \varepsilon_i \nu_i(A)), \qquad \varepsilon_i = \varepsilon_{\gamma_i}$ and $\nu_i = \nu_{\gamma_i}$,

where $\varepsilon_\gamma$ and $\nu_\gamma$ are defined in (34). We may now state the main result of this section.

THEOREM 12. *Assume* (SA0)–(SA2). *For* $\xi \geq 0$, *set*

(42) $$\gamma_i = \frac{\log(i+1)}{d(1+\xi)} + \underline{\gamma},$$

*where $d$ is defined in* (35). *Then for any initial probability measure $\mu$, we have*

(43) $$\lim_{n \to \infty} \|\mu P_1 \cdots P_n - \pi_n\|_{\mathrm{TV}} = 0.$$

PROOF. For any $1 \leq m \leq n$, we have

$$\|\mu P_1 \cdots P_n - \pi_n\|_{\mathrm{TV}}$$
(44)
$$\leq \|(\mu P_1 \cdots P_m)P_{m+1} \cdots P_n - \pi_m P_{m+1} \cdots P_n\|_{\mathrm{TV}}$$
$$+ \sum_{l=m}^{n-1} \|\pi_l P_{l+1} P_{l+2} \cdots P_n - \pi_{l+1} P_{l+1} P_{l+2} \cdots P_n\|_{\mathrm{TV}}.$$

Let $(a_n, n \geq 0)$ be a sequence of integers such that $\limsup_{n \to \infty}(a_n^{-1} + a_n/n) = 0$. Note that for sufficiently large $n$,

$$(\lambda^{\mathrm{Leb}}(C))^{-1} \sum_{i=n-a_n}^{n} \varepsilon_i = \varepsilon \sum_{i=n-a_n}^{n} e^{-\gamma_i d} = e^{-\underline{\gamma} d} \varepsilon \sum_{i=n-a_n}^{n} (1+i)^{-1/(1+\xi)}.$$

Hence $\lim_{n \to \infty} \sum_{i=n-a_n}^{n} \varepsilon_i = \infty$.

From Proposition 11, we have $\sup_i \sup_{(x,x') \in \bar{C}} \bar{R}_i \bar{V}_s(x,x') < \infty$, and thus there exists an integer $l$ such that $\lambda^l \sup_i \sup_{(x,x') \in C} \bar{R}_i \bar{V}_s(x,x') \leq \lambda$, with $\lambda < 1$ satisfying (40). Since

$$\lambda^l \sup_i \sup_{(x,x') \in C} \bar{R}_i \bar{V}_s(x,x') \lambda^{-1} \leq 1,$$

Theorem 8 implies that, for all $n \geq (l+1)a_n$ and any initial distributions $\xi$ and $\xi'$,

$$\|\xi P_{n-(l+1)a_n} \cdots P_n - \xi' P_{n-(l+1)a_n} \cdots P_n\|_{\text{TV}}$$
$$\leq \left[\prod_{i=n-a_n}^{n}(1-\varepsilon_i)\right] + \lambda^{a_n} \xi \otimes \xi'(\bar{V}_s)$$
$$\leq \exp\left(-\sum_{i=n-a_n}^{n} \varepsilon_i\right) + \lambda^{a_n} \xi \otimes \xi'(\bar{V}_s).$$

To bound the first term on the right-hand side of (44), we use the expression above with $\xi = \mu P_1 \cdots P_m$ and $\xi' = \pi_m$ with $m = n - (l+1)a_n - 1$. Equation (39) implies that for any initial distribution $\mu$ and any integer $m$,

$$\mu P_1 \cdots P_m V_s \leq \lambda_0^m \mu V_s + \frac{b}{1-\lambda_0}.$$

Since $\pi_m P_m = \pi_m$,

$$\pi_m V_s \leq \lambda_0 \pi_m V_s + b \quad \Longrightarrow \quad \pi_m V_s \leq \frac{b}{1-\lambda_0}.$$

Hence $\mu P_1 \cdots P_m \otimes \pi_m(\bar{V}_s) \leq \lambda_0^m \mu V_s/2 + b/(1-\lambda_0) < \infty$, which implies

$$(45) \quad \lim_{n \to \infty} \|(\mu P_1 \cdots P_{n-(l+1)a_n-1})P_{n-(l+1)a_n} \cdots P_n - \pi_{n-(l+1)a_n-1} P_{n-(l+1)a_n} \cdots P_n\|_{\text{TV}} = 0.$$

We now bound the second term on the right-hand side of (44). For any $l \in \{1, \ldots, n\}$, $\|\pi_l P_{l+1} \cdots P_n - \pi_{l+1} P_{l+1} \cdots P_n\|_{\text{TV}} \leq \|\pi_l - \pi_{l+1}\|_{\text{TV}}$ and thus

$$\sum_{l=m}^{n-1} \|\pi_l P_{l+1} \cdots P_n - \pi_{l+1} P_{l+1} \cdots P_n\|_{\text{TV}} \leq \sum_{l=m}^{n-1} \|\pi_l - \pi_{l+1}\|_{\text{TV}}.$$

To bound this difference we use Lemma A.1, which simplifies the argument in Haario, Saksman and Tamminen (2001). This lemma shows that

$$(46) \quad \sum_{l=m}^{n-1} \|\pi_l - \pi_{l+1}\|_{\text{TV}} \leq 2\log(Z(\gamma_m)/Z(\gamma_n)),$$

where $Z(\gamma) = \int_{\mathbb{R}} e^{-\gamma f(x)} \mu^{\text{Leb}}(dx) / \sup_{x \in \mathbb{R}} e^{-\gamma f(x)}$. Using the Laplace formula [see, e.g., Barndorff-Nielsen and Cox (1989)], it may be shown that

$$(47) \quad Z(\gamma) = (2\pi\gamma^{-1})^{1/2} \left( \sum_{x \in \mathcal{M}} (f''(x))^{-1/2} \right) (1 + o(1)) \qquad \text{as } \gamma \to \infty,$$

where $\mathcal{M}$ is the set of global minima of $f(x)$ (recall that, under the stated assumptions, these minima are isolated and there are only a finite number of them). For any integer $j$, (46) and (47) show that

$$(48) \quad \lim_{n \to \infty} \sum_{l=n-ja_n}^{n-1} \|\pi_l - \pi_{l+1}\|_{\text{TV}} \leq 2 \lim_{n \to \infty} \log\left(\frac{Z(\gamma_{n-ja_n})}{Z(\gamma_n)}\right)$$
$$\leq \lim_{n \to \infty} \log\left(\frac{\gamma_n}{\gamma_{n-ja_n}}\right) = 0.$$

Together with (45), this concludes the proof. □

## APPENDIX A: TECHNICAL LEMMAS

LEMMA A.1. *Let $h$ be a nonnegative function on a measurable space $(\mathcal{X}, \mathcal{B}(\mathcal{X}), \mu)$. Assume that $0 < \int h^\gamma \, d\mu < \infty$ for all $\gamma > \gamma_0 > 0$ and that $\|h\|_\infty = \operatorname{ess\,sup}_\mathcal{X} h(x) := \inf\{M : \mu\{x : h(x) > M\} = 0\} < \infty$. For $\gamma \geq \gamma_0$, denote by $\mu_\gamma$ the measure over $(\mathcal{X}, \mathcal{B}(\mathcal{X}))$ with probability density function $h^\gamma / \int h^\gamma \, d\mu$ w.r.t. $\mu$. Then, for $\gamma' \geq \gamma \geq \gamma_0$,*

$$\|\mu_\gamma - \mu_{\gamma'}\|_{\text{TV}} \leq 2 \log\left(\frac{Z(\gamma)}{Z(\gamma')}\right), \qquad Z(\gamma) = \frac{\int h^\gamma \, d\mu}{\|h\|_\infty^\gamma}.$$

PROOF. Sheff's identity shows that

$$\|\mu_\gamma - \mu_{\gamma'}\|_{\text{TV}} = \int |f - g| \, d\mu,$$

where $f = h^\gamma / \int h^\gamma \, d\mu$ and $g = h^{\gamma'} / \int h^{\gamma'} \, d\mu$. Note that $f / \|f\|_\infty = (h/\|h\|_\infty)^\gamma \geq g/\|g\|_\infty = (h/\|h\|_\infty)^{\gamma'}$, $\mu$-a.e. and $\|g\|_\infty / \|f\|_\infty = Z(\gamma)/Z(\gamma')$. The proof follows from Lemma A.2, which may be of independent interest. □

LEMMA A.2. *Let $f$ and $g$ be two probability density functions w.r.t. a common dominating measure $\mu$ on $(\mathcal{X}, \mathcal{B}(\mathcal{X}))$. Assume that $\|f\|_\infty < \infty$ and $\|g\|_\infty < \infty$, and $f(x)/\|f\|_\infty \geq g(x)/\|g\|_\infty$, $\mu$-a.s. Then*

$$\int |f - g| \, d\mu \leq 2 \log(\|g\|_\infty / \|f\|_\infty).$$

PROOF. Using the inequality $(\|f\|_\infty/\|g\|_\infty)g \leq f$ and $|f-g| = f + g - 2(f \wedge g)$, we have

$$\int |f - g|\, d\mu = 2\left(1 - \int f \wedge g\, d\mu\right)$$
$$\leq 2\left(1 - \int \frac{\|f\|_\infty g}{\|g\|_\infty} \wedge g\, d\mu\right)$$
$$= 2\left(1 - \frac{\|f\|_\infty}{\|g\|_\infty}\right)$$

and the proof follows from the inequality

$$1 - x \leq \log(1/x) \qquad \text{for } x > 0. \qquad \square$$

**Acknowledgments.** We sincerely thank the referee and an Associate Editor for their very careful reading of the manuscript which led to many improvements.

R. Douc
CMAP
École Polytechnique
Route de Saclay
91128 Palaiseau Cedex
France
e-mail: douc@cmapx.polytechnique.fr

E. Moulines
Département TSI/CNRS URA 820
Ecole Nationale Supérieure
  des Télécommunications (ENST)
46 Rue Barrault
75634 Paris Cedex 13
France
e-mail: moulines@tsi.enst.fr

J. S. Rosenthal
Department of Statistics
University of Toronto
Toronto, Ontario
Canada M5S 3G3
e-mail: jeff@math.toronto.edu